 \documentclass[12pt]{amsart}
\usepackage[latin1]{inputenc}
\usepackage[english]{babel}
\usepackage{amssymb,amsmath,amsthm}
\usepackage{amsfonts}
\usepackage{fancyhdr}

\newtheorem{theorem}{Theorem}[section]
\newtheorem{proposition}{Proposition}[section]

\newtheorem{lemme}{Lemma}[section]

\newtheorem{Remarque}{Remark}[section]

\begin{document}

\title{On the completed tensor product of two algebras on a field}

\author{Mohamed TABA\^A}
\address{DEPARTMENT OF MATHEMATICS\\ FACULTY OF SCIENCES\\
MOHAMMED V UNIVERSITY IN RABAT\\ RABAT\\ MOROCCO
}
\email{mohamedtabaa11@gmail.com; tabaa@fsr.ac.ma}

\begin{abstract}
We give a response to a question posed by Groethendieck on the transfert of the  properties: reduced, normal, domain, regular,  complete intersection,  Gorenstein, Cohen-Macaulay,  to the completed tensor product of two Noetherian algebras on a field.
\end{abstract}

\maketitle

%%%%%%%%%%%%%%%%%%%%%%%%%%%%%%%%%%%%%%%%%%%%%%%%%%%%%%%%%
\footnotesize{\noindent{\bf Key Words.}   Completed tensor product,  reduced,   domain,     complete intersection.  \medskip 

 \noindent{\bf  Mathematics Subject Classification (2010)}  13J10, 13H10, 13G05}
\bigskip\bigskip
 %%%%%%%%%%%%%%%%%%%%%%%%%%%%%%%%%%%%%%

\date{\today}
%\end{frontmatter}

\section{Introduction}
 Throughout this paper, all rings
are commutative with   identity     and all modules are unitary.  \medskip

it is shown in  \cite{6} that, if $k$ is a perfect  (resp.   algebraically closed) field, $A, B$ two complete Noetherian local rings containing $k$ and whose residual fields are   finite extensions of $k$, then, if $A$ and $B$ are reduced (resp. domains), so is the completed tensor product  $ A \hat{\otimes} _{k}B$.\medskip
 
At the end of  \cite[Remark  7.5.8.]{6}, it is written that ``il serait d\'{e}%
sirable de d\'{e}velopper des variantes de (7.5.6) et (7.5.7) o\`{u} on
affaiblirait l'hypoth\`{e}se de finitude des corps r\'{e}siduels de A et B,
en supposant par exemple qu'un seul d'entre eux est extension finie de k,
l'autre \'{e}tant quelconque".  In this paper we prove that,  if $k$ is perfect  and $ A \hat{\otimes} _{k}B$ is  Noetherian, then  $A \hat{\otimes} _{k}B$ is reduced    (resp. normal) if $A$ and $B$ are so, and if     $k$ is algebraically closed  and $ A \hat{\otimes} _{k}B$ is  Noetherian, then 
$A \hat{\otimes} _{k}B$ is a domain if  $A$ and $B$ are so.   \medskip

On the other hand, in    \cite[Theorem 6]{W}, Watanabe, Isikakawa,
Tachibana and Otsuka show that, if $A$ and $B$ are two  Gorenstein 
semilocal rings  containing a field $k$,  then the  completed tensor product  $ A \hat{\otimes} _{k}B$  is a Gorenstein ring if, for any   maximal ideal $\mathfrak{n}$ of $B$,   $ B /  \mathfrak{n}$ is a finite extension  of $k$.  In this paper, we show  that if $ A \hat{\otimes} _{k}B$    is Noetherian, then it is     complete intersection  (resp. Gorenstein,   Cohen-Macaulay) rings  if $A$ and $B$ are so, and if furthermore $k$ is perfect, then  it  is regular
   if    $A$ and $B$ are so.

%%%%%%%%%%%%%%%%%%%%%%%%%%%%%%%%%%%%%%%%%%%%%%%%%%%%%%%%%%%%%%%%%%%%%%%%%%%%%%%%%%%%%%%%%%
%%%%%%%%%%%%%%%%%%%%%%%%%%%%%%%%%%%%%%%%%%%%%%%%%%%%%%%%%%%%%%%%%%%%%%%%%%%%%%%%%%%%%%%%%%
%%%%%%%%%%%%%%%%%%%%%%%%%%%%%%%%%%%%%%%%%%%%%%%%%%%%%%%%%%%%%%%%%%%%%%%%%%%%%%%%%%%%%%%%%% 
\section{Preliminaries}
Recall the two theorems of \cite{10} in the form we are going to use.

\begin{theorem} \label{thm-0.1}
Let $B$ be a local Noetherian rings, $\mathfrak{n}$ its 
maximal ideal and $M$ a $B$-module. If $M$ is flat on $B$, then its  
 separated completion for the $\mathfrak{n}$-adic topology is also flat on $B$.
\end{theorem}
  
\begin{theorem} \label{thm-0.2}
Let $ A $ and $ B $ two local Noetherian rings  containing a field $k$. If the ring $A \hat{\otimes} _{k}B$ is Noetherian, then  it is flat on both  $A$ and   $B$.
\end{theorem}
 
 In what follows, we freely use the following proposition:

\begin{proposition} \label{Prop-0.1}
  Let $ A $ and $ B $ be two rings containing a field $k$,  $\mathfrak{m}$ is the radical of $A$, $\mathfrak{n}$  is the radical of $B$, $C $ the ring $A\otimes_{k}B$, $\mathfrak{r}$ the ideal $\mathfrak{m}\otimes _{k}B+ A\otimes _{k}\mathfrak{n}$ of $A\otimes _{k}B$  and 
  $E$  the ring   $ A \hat{\otimes} _{k}B$.
 Then:
\begin{itemize}
\item[(i)]   $ \widehat{\mathfrak{r}} $ (and a fortiori $\mathfrak{m}E$ and $\mathfrak{n}E$) is contained in the radical of $E$.
\item[(ii)]  If furthermore $\mathfrak{m}$ and $\mathfrak{n}$ are finitely generated, then $ \widehat{\mathfrak{r}}=\mathfrak{r}E $  and the topology of $E$ is the topology
$\mathfrak{r}$-adic.
\end{itemize}
\end{proposition}
\begin{proof}
 We have $\mathfrak{r}^{2i}\subset  \mathfrak{m}\otimes _{k}B+ A\otimes _{k}\mathfrak{n} \subset \mathfrak{r}^{i}$. So $E$ is the   separated completion of $ C$ for the $\mathfrak{r}$-adic topology. We deduce that 
  $\widehat{\mathfrak{r}} $   is contained in the radical of $E$   and  that $ \widehat{\mathfrak{r}^i}=\mathfrak{r}^iE $  if  $\mathfrak{m}$ and $\mathfrak{n}$ are  finitely generated.
\end{proof}

\section{Main results}
 
In what follows we freely use use without explicit mention the results  of  \cite{1} and  \cite{4}.
  The notations we will use in the proofs are those  of Proposition \ref{Prop-0.1}.

\subsection{Case: Complete intersection, Gorenstein,   Cohen-Macaulay, Regular}
 
\begin{theorem} \label{thm-1.1}
   Let $ A $ and $ B $ be  two semilocal Noetherian  rings containing a field $k$ such  that the ring  $ A \hat{\otimes} _{k}B$   is Noetherian. Then:
\begin{itemize}
\item[(i)]  If $A$ and $B$ are     complete intersection      (resp.
Gorenstein, resp. Cohen-Macaulay) rings, so is  $ A \hat{\otimes} _{k}B$.
\item[(ii)]  Suppose further that $k$ is perfect. If $A$ and $B$ are regular, so is   
  $ A \hat{\otimes} _{k}B$.
\end{itemize}
\end{theorem} 
\begin{proof} Case where $A$ and $B$ are local:  We designate by \textit{\textbf{P}} one of the following properties of Noetherian rings:   complete intersection, Gorenstein, Cohen-Macaulay.

 (i) Let $\mathfrak{Q}$ be a  maximal ideal of $E$. From Theorem \ref{thm-0.2}, $E$ is flat on $A$ and, as $ \mathfrak{m}E$ is contained in $\mathfrak{Q}$, the homomorphism $ A \rightarrow
E_ {\mathfrak{Q}} $ is local and flat. To show that $ E_ {\mathfrak{Q}} $  verifies \textit{\textbf{P}}   it is enough
to show that the closed fiber   $ E_ {\mathfrak{Q}} / \mathfrak{m}E_ {\mathfrak{Q}} $  verifies \textit{\textbf{P}}   because $A$   verifies \textit{\textbf{P}}. 
We claim that the  composed homomorphism $B\to E\to  E/  \mathfrak{m} E$, which is by definition the composed homomorphism  $B\to C / \mathfrak{m} C\to  E/  \mathfrak{m} E$, is flat.   As $E$ is a Zariski ring for the $\mathfrak{r}$-adic topology, the ideal $ \mathfrak{m}E$ is closed and  since the topology of $ C/  \mathfrak{m} C $ induced by that of $C$ is the $\mathfrak{n}$-adic topology,  $ E /  \mathfrak{m}E $ is the separated  completion  of $ C / \mathfrak{m} C $ for the $ \mathfrak{n}$-adic topology. But $ C/  \mathfrak{m} C$ is isomorphic to  $(A / \mathfrak{m}) \otimes _ {k} B $, so it is flat on $B$. According to Theorem \ref{thm-0.1}, $ E / \mathfrak{m}E $  is flat on $B$. We deduce that the homomorphism $ B \rightarrow E_ {\mathfrak{Q}} / \mathfrak{m}E_ {\mathfrak{Q}} $ is flat.   Since $\mathfrak{n}E$  is contained in $\mathfrak{Q}$,  it is local.  Its closed fiber is isomorphic to $ (E / \mathfrak{r}E) _ {\mathfrak{Q} / \mathfrak{r}E} $.  But $ E / \mathfrak{r}E $ is isomorphic  to  $
(A / \mathfrak{m}) \otimes _ {k} (B / \mathfrak{n}) $ and according to  Proposition 2.7 of \cite{9}, $ (A / \mathfrak{m}) \otimes _ {k} (B / \mathfrak{n}) $  verifies \textit{\textbf{P}}, so $ E / \mathfrak{r}E $  verifies \textit{\textbf{P}}.  We deduce that  $ (E / \mathfrak{r}E) _ {\mathfrak{Q} / \mathfrak{r}E} $  verifies \textit{\textbf{P}}   and hence $ E_ {\mathfrak{Q}} /  \mathfrak{m}E_ {\mathfrak{Q}} $  verifies \textit{\textbf{P}}   because $B$  verifies \textit{\textbf{P}}.
  
(ii) By Corollary 2.6 of \cite{9},  the ring $ (A / \mathfrak{m}) \otimes _ {k} (B / \mathfrak{n}) $ is   regular. It suffices to  replace in the previous proof \textit{\textbf{P}}  by regular.

General case:
 Let $( \mathfrak{m}_ {r} )_r$   the distinct maximal ideals 
 of  $A$ and  $ (\mathfrak{n}_ {s}  )_s$  the distinct maximal ideals   of  $B$.

We prove that 
$$ A \hat{\otimes} _{k}B \cong \prod\limits_ {r,s} ( A_{\mathfrak{m}_{r} }\hat{\otimes} _{k}B_{ \mathfrak{n}_{s}}).$$
We have
$$A/\mathfrak{m}^{i}=\prod\limits_ {r}  (A/\mathfrak{m}_{r}^{i})=\prod\limits_ {r}
(A_{\mathfrak{m}_{r}}/\mathfrak{m}_{r}^{i}A_{m_{r}})\ \mathrm{and} \ B/ \mathfrak{n}^{j}=\prod\limits_ {s}  (B/ \mathfrak{n}_{s}^{j})=\prod\limits_ {s} 
(B_{\mathfrak{n}_{s}}/ \mathfrak{n}_{s}^{j}B_{ \mathfrak{n}_{s}}).$$
Whence 
\begin{eqnarray*}
A \hat{\otimes} _{k}B &=& \underset{i,j}{\underleftarrow{\lim }}((A/\mathfrak{m}^{i})\otimes _{k} (B/ \mathfrak{n}^{j}))\\
 &=&  \underset{i,j}{\underleftarrow{\lim }}(\prod\limits_ {r}  
(A_{\mathfrak{m}_ {r} }/\mathfrak{m}_{r}^{i}A_{\mathfrak{m}_{r}})\otimes _{k} \prod\limits_ {s}   (B_{ \mathfrak{n}_{s}}/\mathfrak{n}_{s}^{j}B_{ \mathfrak{n}_{s}}))\\
&=&\underset{i,j}{\underleftarrow{\lim }} \,  \prod\limits_ {r,s}
((A_{\mathfrak{m}_{r}}/\mathfrak{m}_{r}^{i}A_{\mathfrak{m}_ {r} })\otimes _{k} (B_{\mathfrak{n}_{s}}/ \mathfrak{n}_{s}^{j}B_{ \mathfrak{n}_{s}}))\\
&=& \prod\limits_ {r,s}  \, \underset{i,j}{\underleftarrow{\lim }} 
((A_{\mathfrak{m}_ {r} }/\mathfrak{m}_{r}^{i}A_{\mathfrak{m}_{r}})\otimes _{k} (B_{\mathfrak{n}_{s}}/ \mathfrak{n}_{s}^{j}B_{ \mathfrak{n}_{s}}))\\
&=& 
\prod\limits_ {r,s}  ( A_{\mathfrak{m}_{r} }\hat{\otimes} _{k}B_{ \mathfrak{n}_{s}}).
\end{eqnarray*} 

For all $r, s$, the ring $ A_{\mathfrak{m}_{r} }\hat{\otimes} _{k}B_{ \mathfrak{n}_{s}}$  is Noetherian, so it suffices to apply the previous case.
\end{proof}

\begin{Remarque} \label{rem.1}  Taking account of     \cite[Err$_{IV}$,\, 19]{66}, the proof of the theorem shows also that $$dim(A \hat{\otimes} _{k}B)=sup_ {r,s}\{ht(\mathfrak{m}_r)+ht(\mathfrak{n}_s)+inf(deg.tr_k(A/\mathfrak{m}_r),deg.tr_k(B/\mathfrak{n}_s))\}.$$
If furthermore, for all $r,s$, one of the extension $A/\mathfrak{m}_r$ or $B/\mathfrak{n}_s$ is algebraic over $k$ then
$$dim(A \hat{\otimes} _{k}B)=dim(A)+dim(B).$$
\end{Remarque}
%%%%%%%%%%%%%%%%%%%%%%%%%%%%%%%%

\subsection{Case: Reduced, Normal.} Recall, \cite[7.3.1.]{6},  that a homomorphism of Noetherian  rings    $\sigma :A\rightarrow B$  is said to be reduced  (resp.
normal) if it is flat and if, for any prime ideal $ \mathfrak{p} $ of $A$, the $ k (\mathfrak{p} ) $-algebra  $B\otimes _{A}k(\mathfrak{p} )$  is geometrically reduced (resp. geometrically normal).

\begin{lemme} \label{lem.1}  Let $A$ and  $B$ be two complete Noetherian local rings  
containing a field  $k$ such that   $A \hat{\otimes} _{k}B$  is Noetherian. Let $\mathfrak{m}$ be the maximal ideal of $A$. For any finitely generated $A$-module  $M$ endowed  with  the $\mathfrak{m}$-adic topology, the completed tensor product  $M \hat{\otimes} _{k}B$ is identified with $M\otimes _{A}(A \hat{\otimes} _{k}B )$.
\end{lemme}

\begin{proof} According to  \cite[Proposition 2.1]{10}  applied 
to the $A$-algebra $E$,  $M\hat{\otimes}_{A}E$  is  identified  with  $%
M\otimes _{A}E$.  It suffices  to apply  \cite[Proposition 7.7.8.]{5}.
\end{proof}

Recall the following localization theorem: 

\begin{theorem}[\cite{7}] \label{lem-2} 
  Let $ \sigma: A \rightarrow B $ be a flat  homomorphism of local Noetherian rings and  $\mathfrak{m}$
the maximal ideal of $A$. We assume that the formal  fibers   of $A$ are geometrically reduced  (resp. geometrically normal). If the closed fiber $B/\mathfrak{m}B$  is geometrically reduced  (resp. geometrically normal) then  $ \sigma$ is reduced  (resp. normal).
\end{theorem}

\begin{lemme} \label{lem-2} 
  Let $ \sigma: A \rightarrow B $ be a flat homomorphism
of Noetherian rings. Suppose $A$ is a complete local ring with a maximal ideal $\mathfrak{m}$ and
that $ \mathfrak{m}B$ is contained in the radical of $ B $. If the  fiber $ B / \mathfrak{m}B $   is geometrically reduced  (resp. geometrically normal) on $ A / \mathfrak{m}$, then $B$ is reduced (resp. normal) if $A$ is so.
\end{lemme}

\begin{proof}
  Let $\mathfrak{Q}$ be a maximal ideal of $B$. We have $\sigma ^{-1}(\mathfrak{Q})=\mathfrak{m}$, 
so the homomorphism  $\widetilde{\sigma }:A\rightarrow B_{\mathfrak{Q}}$   deduced from 
$ \sigma $ is local and flat. Its closed fiber is $(B/\mathfrak{m}B)_{\mathfrak{Q}/\mathfrak{m}B}$. It is geometrically reduced  (resp. geometrically normal) on $ A / \mathfrak{m} $.  Using localization theorem,
the homomorphism $ \widetilde {\sigma} $ is reduced  (resp. normal). So $ B_{\mathfrak{Q}} $ is reduced  (resp. normal).
\end{proof}

 \begin{theorem} \label{thm-main1} 
   Let $k$ be a perfect field, $A$ and $B$  two complete Noetherian local rings   
  containing  $k$ such that the ring $A \hat{\otimes} _{k}B$  is Noetherian. If $A$ and $B$ are reduced    (resp. normal), then  so is   $A \hat{\otimes} _{k}B$.
\end{theorem}
 \begin{proof}
From Theorem  \ref{thm-0.2}, $E$  is flat on $A$. Taking account of  Lemma \ref{lem-2}, to show that $E$  is reduced  (resp. normal), it suffices to show that the fiber $ E / \mathfrak{m}E $ is geometrically reduced (resp. geometrically normal) on $ A / \mathfrak{m}, $ because $\mathfrak{m}E$ is contained  in the radical of $E$ and $A$ is reduced  (resp. normal).  Let $L$  be
a finite extension of $ A / \mathfrak{m} $.  Let us show that  $ L\otimes _ {A / \mathfrak{m}} E / \mathfrak{m}E = L\otimes _ {A } E $ is 
reduced (resp. normal). For this, consider the composed homomorphism $B\rightarrow  E \rightarrow   L \otimes _ {A } E $ and prove that it is flat. The following diagram 
$$\begin{array}{ccccc}
B&\rightarrow  &E    \\
\downarrow  &   & \downarrow   \\
L \otimes _ {k } B&\rightarrow& L \otimes _ {A } E
\end{array}
$$ 
is commutative, it suffices then to show that the composed homomorphism $B\rightarrow  L \otimes _ {k } B \rightarrow   L \otimes _ {A } E $ is flat. From  Lemma \ref{lem.1}, $ L \otimes _ {A } E $ is isomorphic to the completed  tensor product $L\hat{\otimes} _{k} B$.  But  $\mathfrak{m}L=0$, so the tensor product topology of $
L \otimes _ {k} B$ is the $\mathfrak{n}$-adic topology. So $ E \otimes _ {A} L $ is the separated  
completion of $ L \otimes _ {k} B $ for this topology.  But $
L \otimes _ {k} B $ is flat on $B$, so by  Theorem \ref{thm-0.1}, $
L \otimes _{A} E$ is flat on $B$.    Let us show that $ \mathfrak{n} (L \otimes _ {A } E ) $ is contained in the radical of $L \otimes _ {A } E$. Let $\mathfrak{Q}$ be a   maximal ideal of $L \otimes _ {A  } E  $.  Since the homomorphism $ E  \rightarrow L \otimes _ {A } E$ is integral,  the inverse image $\mathfrak{Q}'$ of $\mathfrak{Q}$ is a maximal ideal of $E$. Then, $\mathfrak{Q}'$ contains   $  \mathfrak{n}E $ and so $\mathfrak{Q}$ contains $ \mathfrak{n} (L \otimes _ {A } E ) $. It remains to show that the fiber $L  \otimes _ {A } E  \otimes _ {B} (B/ \mathfrak{n})$ is geometrically reduced   (resp. geometrically normal) on $ B / \mathfrak{n} $. By Lemma \ref{lem.1}, the fiber  $L  \otimes _ {A } E  \otimes _ {B} (B/ \mathfrak{n})$ is isomorphic to  $L  \otimes _ {k }   (B/ \mathfrak{n})$, because the  $\mathfrak{m}$-adic topology on $L$ is the discrete topology and the he $\mathfrak{n}$-adic topology on $B/ \mathfrak{n}$ is also  the discrete topology, and since  the field $ k $ is perfect, it is geometrically  regular over $B/ \mathfrak{n}$.    
    It suffices to use Lemma  \ref{lem-2} since $B$ is reduced  (resp. normal).
\end{proof}

The previous result does not hold true if $k$ is not perfect. Let  $p$ its characteristic, $a$  an element of $k-k^p$, and $K$ the extension   $k[X]/(X^p-a)$ of $k$, where $X$ is an indeterminate over $k$. Take $A=K[[Y_1,...,Y_m]]$ and $B=K[[Y_1,...,Y_n]]$, then the Noetherian ring $A \hat{\otimes} _{k}B$ is not reduced.

%%%%%%%%%%%%%%%%%%%%%%%%%%%%%%%%

\subsection{Case: Domain.}  
 
Recall that if $A$ is a ring, $Spec(A)$ is connected if only if $A$ contains 
no idempotents other than $0$ and $1$.

 \begin{proposition}
Let $k$ be a perfect field, $A$ and $B$ be two complete local
Noetherian  rings containing $k$ such that $A\widehat{\otimes }_{k}B$ is Noetherian, $\mathfrak{m}$ the maximal ideal of $A$ and $\mathfrak{n}$ the maximal ideal of $
B $. We suppose that $A$ and $B$ are integrally closed domains. Then $A
\widehat{\otimes }_{k}B$ is a domain if and only if $Spec((A/\mathfrak{m})\otimes
_{k}(B/\mathfrak{n}))$ is connected.
\end{proposition} 

By Ttheorem \ref{thm-main1}, $E$ is locally a domain. Proposition 10.8 of \cite{8}  
applied to $E$ shows that $E$ is a domain if and only if $Spec(E/\mathfrak{r}E)$ is
connected. The equivalence follows from the fact that $E/\mathfrak{r}E$ is isomorphic
to $(A/\mathfrak{m})\otimes _{k}(B/\mathfrak{n})$.

\begin{lemme} Let $A$ and $B$ be two Noetherian local rings
containing a field $k$,  $\mathfrak{m}$ the maximal ideal of $A$ and $\mathfrak{n}$ the maximal
ideal of B. Then $A\widehat{\otimes }_{k}B$ is Noetherian if and only if $%
(A/\mathfrak{m})\otimes _{k}(B/\mathfrak{n})$ is\ Noetherian.\end{lemme} 
\begin{proof}
The ring $E/\mathfrak{r}E$ is isomorphic to $(A/\mathfrak{m})\otimes _{k}(B/\mathfrak{n}).$ So
the sufficient condition is deduced from the fact that the ideal $\mathfrak{r}E$ is  
finitely generated and $E$ is complete.
\end{proof}  

  \begin{theorem}
  Let $k$ be an algebraically closed field, $A$ and 
$B$ be two complete Noetherian local rings containing $k$ such that $A 
\widehat{\otimes }_{k}B$\ is Noetherian. If $A$ and $B$ are domains, then so
is $A\widehat{\otimes }_{k}B.$
  \end{theorem}\begin{proof} 
   By Nagata's theorem (\cite{G6}, 23.1.5 and 23.1.6), the
integral closure $A\prime $ (resp. $B\prime $) of $A$ (resp. $B$) is a
finitely generated $A$-module (resp. $B$-module) and a complete Noetherian
local ring. We show at first, using the previous proposition, that $
A^{\prime }\widehat{\otimes }_{k}B^{\prime }$ is a domain. The ring $A\prime 
$ (resp. $B\prime $) dominates $A$ (resp. $B$) then, if $\mathfrak{m}$ (resp. $\mathfrak{n}$) is the maximal ideal of $A$ (resp. $B$) and $\mathfrak{m}\prime $ (resp. $\mathfrak{n}\prime $) is the maximal ideal of $A\prime $ (resp. $B\prime $),  the homomorphism $
A/\mathfrak{m}\rightarrow A^{\prime }/\mathfrak{m}^{\prime }$ (resp. $B/\mathfrak{n}\rightarrow B^{\prime}/\mathfrak{n}^{\prime }$) induced by the injection $A\rightarrow A^{\prime }$ (resp. $ 
B\rightarrow B^{\prime }$) is finite. We induce that the homomorphism $%
(A/\mathfrak{m})\otimes _{k}(B/\mathfrak{n})\rightarrow (A^{\prime }/\mathfrak{m}^{\prime })\otimes
_{k}(B^{\prime }/\mathfrak{n}^{\prime })$ is also finite. By previous lemma the ring $A^{\prime }\widehat{\otimes }_{k}B^{\prime }$ is Noetherian. It remains to
prove that $Spec((A/\mathfrak{m})\otimes _{k}(B/\mathfrak{n}))$ is connected and this follows from
the fact that the ring $(A^{\prime }/\mathfrak{m}^{\prime })\otimes _{k}(B^{\prime
}/\mathfrak{n}^{\prime })$ is a domain since $k$ is algebraically closed. We show that $%
A\widehat{\otimes }_{k}B$ is identified with a subring of $A^{\prime }%
\widehat{\otimes }_{k}B^{\prime }.$ From Theorem 2.2, $A\widehat{\otimes }%
_{k}B$ is flat over A, then $A\widehat{\otimes }_{k}B$ is   identified with
a subring of $A^{\prime }\otimes _{A}(A\widehat{\otimes }_{k}B)$. But, using Lemma \ref{lem.1}, $A^{\prime }\otimes _{A}(A\widehat{\otimes }_{k}B)$ is isomorphic
to $A^{\prime }\widehat{\otimes }_{k}B$. Then, $A\widehat{\otimes }_{k}B$ is
identified with a subring of $A^{\prime }\widehat{\otimes }_{k}B$. The same
argument shows also that $A^{\prime }\widehat{\otimes }_{k}B$ is identified
with a subring of $A^{\prime }\widehat{\otimes }_{k}B^{\prime }$. We
conclude that $A\widehat{\otimes }_{k}B$ is a domain.
\end{proof}  

\noindent \textbf{Acknowledgment.} I would like to thank Driss Bennis for the help and discussion.

\end{document}